\mag=\magstep1
\documentstyle{amsppt} \input amsppt1

\pageheight{23true cm}
\pagewidth{16true cm}
\parindent=4mm
\parskip=3pt plus1pt minus.5pt
\nologo\NoRunningHeads\NoBlackBoxes

\def\i{\looparrowright}
\def\e{\hookrightarrow}

\def\f{\flushpar }
\def\nl{\newline }
\def\np{\newpage }
\def\x{\times }

\topmatter
\title
The intersection of three spheres in a sphere and 
a new application of the Sato-Levine invariant 
\endtitle
\author
Eiji Ogasa
\endauthor

\thanks{
\nl
This research was partially supported by Research Fellowships 
of the Promotion of Science for Young Scientists.}

\endtopmatter

\document
\baselineskip11pt

\vskip5mm
This paper is published in \newline
 Proceedings of the American Mathematical society  126, 1998, PP.3109-3116.   \newline
This manuscript is not the published version.

\f{\bf Abstract.}   
Take transverse immersions $f:S^4_1\amalg$ $S^4_2\amalg$ $S^4_3\i S^6$ 
such that 
(1) $f\vert S^4_i$ is an embedding, 
(2)
$f(S^4_i)\cap f(S^4_j)$$\neq$$\phi$ 
and  
$f(S^4_i)\cap f(S^4_j)$  is connected, 
 and                                                  
(3)$f(S^4_1)\cap f(S^4_2)\cap f(S^4_3)$ $=\phi$.      
Then we obtain  three surface-links                   
$L_i$=                                                
($f^{-1}(f(S^4_i)\cap f(S^4_j))$,  $f^{-1}(f(S^4_i)\cap f(S^4_k))$ )
in $S^4_i$,                                           
where  $(i,j,k)$=(1,2,3), (2,3,1), (3,1,2).           
We prove that,                                        
we have the equality                                  
$\beta(L_1)+$ $\beta(L_2)+$ $\beta(L_3)=0,$           
where                                                 
$\beta(L_i)$ is the Sato-Levine invariant of $L_i$,   
if all $L_i$ are semi-boundary links.

\head 1. Introduction and Main results \endhead

Take transverse immersions $f:S^4_1\amalg$ $S^4_2\amalg$ $S^4_3\i S^6$ 
such that 
(1) 
    $f\vert S^4_i$ is an embedding, 
(2)
$f(S^4_i)\cap f(S^4_j)$$\neq$$\phi$ 
and  
$f(S^4_i)\cap f(S^4_j)$  is connected, 
 and 
(3)$f(S^4_1)\cap f(S^4_2)\cap f(S^4_3)$ $=\phi$. 
Then we obtain  three surface-links 
$L_i$=
($f^{-1}(f(S^4_i)\cap f(S^4_j))$,  $f^{-1}(f(S^4_i)\cap f(S^4_k))$ )
in $S^4_i$, 
where  $(i,j,k)$=(1,2,3), (2,3,1), (3,1,2). 
An orientation is given to each naturally.   
In this paper, we discuss which ones we obtain.

In order to state our theorems, we need some definitions.

We work in the smooth category. 
$S^4_i$ $\cap$ $S^4_j$ is a closed orientable connected surface and 
is oriented naturally. 
Hereafter, 
a {\it surface} will always mean a closed oriented connected surface unless otherwise stated. 

A {\it surface-$( F_1,...,F_\mu)$-link } is a submanifold 
$L=(K_1,...,K_\mu)$ of $S^4$ 
such that $K_i$ is diffeomorphic to the oriented surface $F_i$. 
If $\mu=1$, $L$ is called a {\it surface-$F_1$-knot}.  
A surface-$( F_1, F_2 )$-link  $L=(K_1, K_2)$ is called 
a {\it semi-boundary link} if 
$[K_i]=0$ $\in$ $H_2(S^4-K_j;\Bbb Z)$ $(i\neq j)$ 
 ([18]).   
A surface-$( F_1, F_2 )$-link $L=(K_1, K_2)$ 
is called a {\it boundary link} if 
there exist Seifert hypersurfaces $V_i$ for $K_i$ ($i=1,2$) 
such that $V_1\cap V_2$=$\phi$.
A surface-($F_1, F_2$)-link $(K_1, K_2)$ is called a {\it split link} 
if there exist 4-balls $B^4_1$ and $B^4_2$ in $S^4$ such that 
 $B^4_1\cap$$B^4_2=\phi$ and $K_i\subset B^4_i$.

\definition{Definition}
$(L_1, L_2,L_3)$ is called a {\it triple of surface-links} 
if  
$L_1$ is a $(F_{12},F_{13})$-link, 
$L_2$ is a $(F_{23},F_{21})$-link, 
$L_3$ is a $(F_{31},F_{32})$-link, 
and 
$F_{ij}$ is diffeomorphic to $F_{ji}$ 
($(i,j)$=(1,2),(2,3),(3,1)).
\enddefinition

\definition{Definition}
A triple of surface-links $(L_1, L_2,L_3)$ is said to be {\it realizable} 
if there exists a transverese immersion 
$f:S^4_1\amalg S^4_2\amalg S^4_3\i S^6$ 
such that 
(1)
   $f\vert S^4_i$ is an embedding (i=1,2,3), 
(2)
   ($f^{-1}(f(S^4_i)\cap f(S^4_j))$, $f^{-1}(f(S^4_i)\cap f(S^4_k))$ )
   in $S^4_i$  is $L_i$=$(K_{ij},K_{ik})$  
   ( $(i,j,k)$=(1,2,3), (2,3,1), (3,1,2)), 
and  
(3)$f(S^4_1)\cap f(S^4_2)\cap f(S^4_3)$ $=\phi$. 
\enddefinition

We state the main theorem.
\proclaim{Theorem 1.1} 
Let $L_1$, $L_2$ and  $L_3$ be semi-boundary surface-links.  
Let ($L_1$,$L_2$,$L_3$) be a triple of surface-links.  
Suppose the triple of surface-links $(L_1, L_2,L_3)$ is  realizable. 
Then we have the equality 
$$\beta(L_1)+\beta(L_2)+\beta(L_3)=0,$$ 
where $\beta(L_i)$ is the Sato-Levine invariant of $L_i$. 
\endproclaim

We review the Sato-Levine invariants in \S 2.  
Since there exists a triple of surface-links ($L_1$, $L_2$, $L_3$) 
such that  $\beta(L_1)$=$\beta(L_2)$=0 and $\beta(L_3)$=1 
(See  \S 2.), 
we have:
\proclaim{Corollary 1.2} 
Not all triple of surface-links are realizable. 
\endproclaim

We prove: 

\proclaim{Theorem 1.3}  
There exists a realizable triple of surface-links  $(L_1, L_2,L_3)$ 
such that  $\beta(L_1)$=1,  $\beta(L_2)$=1, and  $\beta(L_3)$=0.  
\endproclaim

We prove the following sufficient conditions for the realization. 

\proclaim{Theorem 1.4}  
Let $L_1$, $L_2$ and $L_3$ be split surface-links. 
Let ($L_1$,  $L_2$, $L_3$) be a triple of surface-links. 
Then the triple of surface-links  $(L_1, L_2,L_3)$ is  realizable. 
\endproclaim

\proclaim{Theorem 1.5}  
Suppose $L_i$ are $(S^2, S^2)$-links and $L_i$ are slice links($i=1,2,3$).  
Then the triple of surface-links  $(L_1, L_2, L_3)$ is realizable.
\endproclaim

We give problems. 
\proclaim{Problem 1.6}
(1) Determine the realizable triple of surface-links. 

(2)    
Is the inverse of Theorem 1.1 valid?   

(3)       
Let $L_1$,  $L_2$ and $L_3$ be $(S^2, S^2)$-links. 
Then is   the triple of surface-links  $(L_1, L_2,L_3)$  realizable?
\endproclaim

\f{\bf Note.}
(i)
Using a result of 
[15] 
(See \S 2.), one can show Problem 1.6.(3) follows from  Problem 1.6.(2).        

(ii)
By Theorem 1.5, if the answer to Problem 1.6.(3) is negative, 
then the answer to an outstanding problem: 
``Is every $(S^2, S^2)$-link slice?'' is negative.  
(Refer to  
[5], [6], and [12]   
for the slice problem.)

This paper is organized as follows. 
In \S 2 we review the Sato-Levine invariant. 
In \S 3 we prove Theorem 1.1. 
In \S 4 we prove Theorem 1.3.
In \S 5 we prove Theorem 1.4.
In \S 6 we prove Theorem 1.5.

\head 2. 
The Sato-Levine invariant and spin cobordism
\endhead

The Sato-Levine invariant is defined by Sato 
 (in [18]) 
and Levine (unpublished) 
independently. 
It is easy to prove that the following definition is equivalent to 
theirs.

\definition{Definition }  
Let $L=(K_1, K_2)$ be a semi-boundary surface-$( F_1, F_2 )$-link. 
Then there exist Seifert hypersurfaces $V_i$ for $K_i$ ($i=1,2$) 
such that $V_i\cap K_j=\phi (i\neq j).$ 
Let $v_i$ be the oriented normal bundle of $V_i$ in $S^4$. 
Let $F$ be the oriented closed surface $V_1$ $\cap$ $V_2$.   
$F$ need not be connected.   
Then the congruence 
 $TS^4\vert_F \cong$ $TF\oplus$ $v_1\vert_ F \oplus$ $v_2\vert_ F$ 
induces a spin structure $\sigma$ on $F$. 
We define 
 the {\it Sato-Levine invariant} $\beta(L)$ of $L$  so that 
$\beta(L)$=$[(F,\sigma)]$ $\in$ $\Omega_2^{\roman{spin}}$ $\cong$ $\Bbb Z_2$  for $L$. 
We call $(F,\sigma)$ a {\it special surface for $L$}.
\enddefinition

By 
 [17] and [18]  
 the following holds. 

\proclaim{Theorem }   
([17] and [18])  
Let $F_1$ be an oriented closed connected surface 
not diffeomorphic to the 2-sphere. 
Let $F_2$ be an arbitrary oriented closed connected surface. 
Then there exists a semi-boundary $( F_1, F_2 )$-link 
whose Sato-Levine invariant is one.
\endproclaim

In 
[15]  
Orr proved the following. 
\proclaim{Theorem } 
([15])  
The Sato-Levine invariant of an arbitrary $(S^2, S^2)$-link is zero.
\endproclaim

The Sato-Levine invariant and its generalization are  studied in 
[1],  
[2],  
[3],  
[4],   
[7],  
[8], 
[10], 
[11], 
[16],  
[19],  
[20],  
P.103 of [21],   
 etc. 
[2]  
says that 
the Sato-Levine invariant is connected with 
[9]. 

\head 3
The proof of Theorem 1.1. 
\endhead
Let $L_1=(K_{12},K_{13})$, $L_2=(K_{23},K_{21})$, and $L_3=(K_{31},K_{32})$. 
Let ${f}:S^4_1\coprod S^4_2\coprod S^4_3 \looparrowright S^6$ 
be an immersion  to realize $(L_1, L_2, L_3)$. 
We abbreviate  $f(S^4_i)$ to  $S^4_i$.  
We first prove: 

\proclaim{ Claim }       
  There exist Seifert hypersurfaces $A_i$ for $S^4_i$  ($i=1, 2, 3 $) 
  such that 
$A_1\cap$ $S^4_2\cap$ $S^4_3$ =$\phi$,   
$A_2\cap$ $S^4_3\cap$ $S^4_1$ =$\phi$, and 
$A_3\cap$ $S^4_1\cap$ $S^4_2$ =$\phi$. 
\endproclaim

\f{\bf Proof.}   
Let $S^4_2$$\times$ $D^2$ be a tubular neighborhood of $S^4_2$  in $S^6$. 
Put 
$D^2$=$\{(x,y) \vert$ $x^2+y^2\leqq0 \}$.  
Then $S^4_2$=$S^4_2\x\{(0,0)\}$.
Put 
$I$=$\{(x,y) \vert$ $0\leqq x\leqq1,$ $y=0\}$.  
We can regard 
 $S^4_2$$\times$ $D^2$ 
as the result of rotating $S^4_2\x I$ around the axis $S^4_2$.  

Put $M=$ $(S^4_2\x I)$ $\cap S^4_1$.
As we rotate $S^4_2\x I$ as above, we rotate $M$ as well.  
The result is $(S^4_2\x D^2)$$\cap S^4_1$.

Take  a Seifert hypersurface $A'_1$ for $S^4_1$. 
Then $A'_1$ $\cap$ $S^4_2$ in $S^4_2$ is a Seifert hypersurface $V'_{21}$  for $K_{21}$.   
We can suppose that 
$A'_1$ $\cap$ $(S^4_2\x p)$ in $S^4_2\x p$ is 
the submanifold $V'_{21}$  for each $p\in D^2$.

Since 
$L_2=(K_{23}, K_{21})$ is a semi-boundary link, 
there is 
a Seifert hypersurface  $V_{21}$  for $K_{21}$  
such that  $V_{21}$ $\cap$ $K_{23}$=$\phi$.
Then there exists a compact oriented 4-manifold $W$ in $S^4_2\x I$ 
with the following properties. 
\roster
\item
$W\cap S^4_2$ in $S^4_2$ is the submanifold $V_{21}$.
\item
$W\cap $ $(S^4_2\x\{(1,0)\})$ in $(S^4_2\x\{(1,0)\})$ is the submanifold $V'_{21}$.
\item
$\overline{(\partial W)-V^2_{21}-V^{'2}_{21}}$ is $M$. 
\endroster

When we rotate $S^4_2\x I$ as above, we rotate $W$ together.  
Let $P$ denote what is made from $W$. 
Note that $\partial P$
=$\partial (\overline{A'_1\cap (S^4_2\x D^2)})$
=$\partial (\overline{A'_1-(A'_1\cap (S^4_2\x D^2))})$ $-S^4_1.$
Let 
$A_1$=$\overline{A'_1-A'_1\cap (S^4_2\x D^2)}$ $\cup P$.

Then $A_1\cap$ $S^4_2\cap$ $S^4_3$ =$V_{21}\cap$ $K_{23}$ =$\phi$. 
Note that, when we modify $A'_1$ to obtain $A_1$, 
we don't change $f$. 

Replace (1,2,3) with (2,3,1) (resp. (3,1,2)) in the above proof. 
Then  we obtain $A_2$ (resp. $A_3$). 
We now obtain  $A_1$, $A_2$ and  $A_3$ so that we keep the immersion $f$. 
This completes the proof.
\qed

Put $X=A_1\cap A_2\cap A_2$.  
Put $F_i$=$(\partial X)\cap S^4_i$. 
Then $\partial X$=$F_1\amalg$ $F_2\amalg$ $F_3$.
By using $A_1$, $A_2$, $A_3$ and $S^6$, 
we give $F_i$ (resp. $W$) a spin structure $\sigma_i$ (resp. $\tau$).
Of course $\partial(X,\tau)=\amalg^3_{i=1}(F_i,\sigma_i)$. 
Then ($F_i$, $\sigma_i$) is a special surface for $L_i$.
Therefore, 
$\Sigma_i\beta(L_i)$
= 
$\Sigma_i[(F_i,\sigma_i)]$=
$[\partial(X,\tau)]$=0 $\in\Omega^{\roman{Spin}}_2$.

\head 4
The proof of Theorem 1.3
\endhead

Let $L_1=(K_{12},K_{13})$ 
be the $(T^2, S^2)$-link in [17].  
Let $L_2=(K_{23},K_{21})$ be the $(S^2,T^2)$-link 
obtained by changing the order of $L_1$.   
Let $L_3=(K_{31},K_{32})$ be the trivial $(S^2,S^2)$-link. 
Note  $\beta(L_1)$=$\beta(L_2)$=1 and  $\beta(L_3)$=0.  
It suffices to prove that 
the triple of surface-links  $(L_1, L_2, L_3)$ is realizable.

$(K_1,K_2)$ is called a {\it pair of surface-$F$-knots} if 
both $K_1$ and $K_2$ are  $F$-knots. 
A pair of $F$-knots $(K_1,K_2)$  is said to be {\it realizable} 
if there exists a transverse immersion 
$f:S^4_1\amalg S^4_2$ $\i S^6$ 
such that (1)$f\vert S^4_i$ is an embedding ($i=1,2$), and  (2)
$f^{-1}(f(S^4_1)\cap f(S^4_2))$ in $S^4_i$ is $K_i$($i=1,2$).

We prove: 

\proclaim{Proposition 4.1}  
Let $K$ be a surface-knot. 
Then the pair of surface-knots $(K, K)$ is realizable. 
\endproclaim

Take an embedding $f:S^{4}_1\coprod S^{4}_2 \e S^{6}$.  
There exists a chart $U$ of $S^{6}$ such that

(1) $\phi:U\cong$  
$ \Bbb R^{4}$ $\times$ $\{(u,v)\vert u,v\in  \Bbb R\}$ $\cong$ 
$ \Bbb R^{4}$ $\times \Bbb R_u \times \Bbb R_v$, and 

(2)$U\cap f(S^{4}_1)=\Bbb R^{4}\times\{(u,v)\vert u=0, v=0\}$. Call it $\Bbb R^4_1$. 

\hskip5mm $U\cap f(S^{4}_2) =\Bbb R^{4}\times\{(u,v)\vert u=1, v=0\}$.   Call it $\Bbb R^4_2$.

We prove Lemma 4.2. Obviously it induces Proposition 4.1. 
\proclaim {Lemma 4.2}  
There exists a transverse immersion   
$g:S^{4}_1\coprod S^{4}_2 \looparrowright S^{6}$     
to realize the pair of surface-knots $(K, K)$ 
with the following properties. 
\roster
\item
 $g\vert S^{4}_2=f\vert S^{4}_2$.            
\item
 $g(S^{4}_2)\cap\Bbb R^4\x\{u=0\}\x\Bbb R_v$
=$g(S^{4}_2)\cap\Bbb R^4\x\{u=0\}\x\{v=0\}$.  
\item 
 $g\vert S^{4}_1$ is isotopic to $f\vert S^{4}_1$.            
\endroster            
\endproclaim 

We modify the embedding $f$ to obtain an immersion $g$.

Take any Seifert hypersurface  $V$ for $K$ in $\Bbb R^{4}_1$.  
Let  
$N(V)=V\times$ $\{t\vert$$-1\leqq t\leqq1 \}$ 
be a tubular neighborhood of $V$ in $\Bbb R^4_1$. 
We define the subset $E$ of $N(V) \times \Bbb R_u \times \Bbb R_v$  

\f$=\{(p,t,u,v)\vert p\in V, -1\leqq t\leqq 1 , u\in \Bbb R , v\in \Bbb R \}$ 
so that 

\hskip3mm$E=\{(p,t,u,v)\vert$  
$p\in V,\quad 0\leqq u\leqq\frac{\pi}{2}, \quad 
t=k\cdot \roman{cos}u,\quad 
v=k\cdot \roman{sin}v, \quad-1\leqq k\leqq 1\}$.

Put $X=$$\overline{(\partial E)-N(V)}$ 
and 
$Y=$ $\overline{f(S^{4}_1)-N(V)}$.  
Then $\partial X$=$\partial Y$= $\partial N(V)$.
Put $\Sigma$= $X\cup Y$. 
Then $\Sigma$ is an embedded 4-sphere.  
We define $g\vert S^{4}_1$ so that $g(S^{4}_1)$=$\Sigma$. 
This completes the proof of Lemma 4.2 and therefore Proposition 4.1.

\f{\bf Note.}  
See Figure 4.1 We draw a lower dimensional analogue. 
There, we replace 
$ \Bbb R^{4}$ $\times \Bbb R_u \times \Bbb R_v$ with 
$\Bbb R^{2}$ $\times \Bbb R_u \times \Bbb R_v$.

\vskip3cm
Figure 4.1.   

You can obtain this figure 
by clicking `PostScript' in the right side of 
the cite of the abstract of this paper in arXiv 
(https://arxiv.org/abs/the number of this paper). 

You can also obtain it from the author's website,  
which can be found by typing his name in search engine.

\vskip3cm

By the definition of $L_i$, 
the $T^2$-knots $K_{12}$ and $K_{21}$ are equivalent. 
Therefore there is an immersion  
$g:S^{4}_1\coprod S^{4}_2 \looparrowright S^{6}$     
to realize the pair of $T^2$-knots ($K_{12}$, $K_{21}$).

We prove the following Lemma 4.3. 
Obviously Lemma 4.3 induce Theorem 1.3. 

\proclaim{Lemma 4.3}
There exists a transverse immersion 
$h:S^{4}_1\coprod S^{4}_2 \coprod S^{4}_3 \looparrowright S^{6}$     
to realize $(L_1, L_2, L_3)$ with the following properties. 
\roster
\item
$h\vert_{S^{4}_1\coprod S^{4}_2}=g$
\item
$h(S^4_3)\subset U$.  
$h(S^4_3)$ is the trivial 3-knot.  
\endroster
\endproclaim

\f{\bf Proof.}   
We modify the immersion $g$ to obtain an immersion $g$. 

  Take  $K_{13}$ (resp. $K_{23}$) 
in $\Bbb R^{4}_1$ (resp. $\Bbb R^{4}_2$). 
There is a Seifert hypersurface $V_{12}$ for $K_{12}$ so that $V_{12}\cap K_{13}=\phi$. 
Take $V_{12}$ as a Seifert hypersurface used in the proof of Lemma 4.2.
Recall  $V_{12}$ and  $K_{13}$ are in $\Bbb R^4_1$.

Recall $K_{13}$ and $K_{23}$ are the trivial $S^2$-knots.
Take a 3-ball $B^3_{13}$ (resp. $B^3_{23}$ )
which bounds $K_{13}$ (resp. $K_{23}$ )
in $\Bbb R^{4}_1$ (resp. $\Bbb R^{4}_2$). 
Note that $B^3_{13}$ does not include in $g(S^4_1)$.

Take the 5-ball $B^5$= 
$\{(q,u,v)\vert q\in  B^3$, $-1\leqq u\leqq 2, -2\leqq v\leqq 2\}$ in $U$. 
Suppose $B^5\cap \Bbb R^4_1$=$B^3_{13}$ and 
$B^5\cap \Bbb R^4_2$=$B^3_{23}$.  
Then $(\partial B^5)\cap$ $S^4_1\cap S^4_2$=$\phi$. 
 
Define $h\vert S^4_3$ so that $h(S^4_3)$=$\partial B^5$.  

This completes the proof of Lemma 4.3 and hence Theorem 1.3.

\head 5
The proof of Theorem 1.4 and a relation between  
knot cobordism and the realization of pair of knots 
\endhead

Surface-$F$-knots  $K_0$ and $K_1$ are said to be 
  {\it cobordant} or   {\it concordant} if 
there is a smooth submanifold $W$ of $S^4\x [0,1]$, 
which meets the boundary transversely in $\partial W$, 
is diffeomorphic to $F\x [0,1]$ and meets $S^4\x \{i\}$ in $K_i$ ($i=0,1$).

We prove the following although it may be folklore. 
\proclaim{Theorem 5.1}
Let $F$ be  a closed connected oriented surface. 
Then arbitrary $F$-knots  $K_0$ and $K_1$ are cobordant. 
\endproclaim

\f{\bf Proof. }      
Let $L$ be a split surface-link with components $K_0$ and $-K_1$. 
It suffices to prove: 

\proclaim{Claim} 
There exists a submanifold 
of $S^4$  which is diffeomorphic to 
$F\x$[0,1] 
such that 
$F\x$[0,1] intersects with $\partial B^5$ transversely,  
$F\x$[0,1] $\cap$ $\partial B^5$ =
$F\x\{0\}$ $\amalg$ $F\x\{1\}$, 
and 
($F\x\{0\}$, $F\x\{1\}$) in  $S^4$=$\partial B^5$ is $L$. 
\endproclaim

Let $V$ be a connected Seifert hypersurface for $L$. 
A spin structure on $V$ is induced from the unique one on $S^4$. 
A spin structure on $\partial V$ is induced from the one on $V$. 
Make a closed spin 3-manifold 
$W=V\cup (F\x [0,1])$ 
so that 
the spin structure on $V$ extend to the one on $W$. 
Note $W$ is not a submanifold of $S^4$. 
Since $\Omega_3^{\roman{spin}}=0$,  
there exists a spin 4-manifold $X$ which $W$ spin-bounds.  
Since $V$ and $F\x [0,1]$ are connected, 
we can take a handle decomposition 
$X=$
$(V\x [0,1])$$\cup$(4-dimensional 2-handles $h^2$)
$\cup \{(F\x [0,1])\x [0,1]\}$.
Take 
$V$$\x$ [0,1] in  
$S^4$ $\x$ [0,1]
so that  
$V$ $\x$ $\{t\}$ is in 
$S^4$ $\x$ $\{t\}$. 
Attach the handles $h^2$ to 
$V$ $\x$ $\{1\}$ $\subset$ $S^4$ $\x$ $\{1\}$. 
Then we can attach the 5-dimensional 2-handles 
$\bar{h}^2$=$h^2\x [-1,1]$ to $S^4$ $\x$$\{1\}$ naturally. 
Let 
$Y$=
 $S^4$ $\x$[0,1] $\cup$(the 5-dimensional 2-handles $\bar{h}^2$).
Since the attaching maps of $\bar{h}^2$ are spin preserving diffeomorphisms, 
$Y$ 
is diffeomorphic to 
$(\natural^* S^2\x B^3)-$(the 5-ball). 
$\partial Y$ is a disjoint union  of 
the standard 4-sphere $S^4_0$ 
and $(\sharp^*$$S^2\x S^2$).   
Hence $Y$ is embedded in $B^5$   
so that $S^4_0$ coincides with $\partial B^5$.

Therefore $F\x [0,1]$ $\subset W$ $\subset$ $B^5$ 
and 
the submanifold 
$F\x [0,1]$ satisfies the condition 
in the Claim. 
This completes the proof. 
\qed

It is easy to prove that Theorem 1.4 is equivalent to the following Theorem 5.2. 
We prove:
\proclaim{Theorem 5.2}   
Let $F$ be a closed connected oriented surface.
If $F$-knots  $K$ and $K'$ are cobordant, 
the pair of $F$-knots  $(K, K')$ is realizable. 
\endproclaim

\f{\bf Proof. }     
By Proposition 4.1, 
the pair of $F$-knots $(K, K')$ is realizable. 
Hence it suffices to prove: 

\proclaim{Claim}  
Suppose that a pair of $F$-knots  $(K_1, K_2)$ is realizable. 
Suppose that $K_2$ is cobordant to $K_3$.  
Then $(K_1, K_3)$ is realizable
\endproclaim

 Proof.    
Let ${f}:S^{4}_1\coprod S^{4}_2 \looparrowright S^{6}$
be  an immersion  to realize $(K_1, K_2)$.
We construct an immersion $\widetilde{f}$ 
$:S^{4}_1\coprod S^{4}_2 \looparrowright S^{6}$
to realize   ($K_1$, $K_3$) as follows.   
Put $\widetilde{f}\vert S^4_2=f\vert S^4_2$.

Let $f(S^4_2)$ $\times$ $D^2$ be a tubular neighborhood of $S^4_2$  in $S^6$. 
Put 
$D^2$=$\{(x,y) \vert$ $x^2+y^2\leqq0 \}$.  
Then $f(S^4_2)$=$S^4_2\x\{(0,0)\}$.
Put 
$I$=$\{(x,y) \vert$ $0\leqq x\leqq1,$ $y=0\}$.  
We can regard 
 $f(S^4_2)$$\times$ $D^2$ 
as what is obtained by rotating $f(S^4_2)\x I$ around $f(S^4_2)$ as the axis.

Put $M=$ $(f(S^4_2)\x I)$ $\cap f(S^4_1)$.
We can regard 
 $(f(S^4_2)\x D^2)$ $\cap f(S^4_1)$ as what is made from $M$ as follows: 
When we rotate $(f(S^4_2)\x I)$ as above, we rotate $M$ together.  
What is made from $M$ is $(f(S^4_2)\x D^2)$$\cap$$f(S^4_1)$.

We can suppose that 
$\{f(S^4_2\x p)\}$ $\cap$ $f(S^4_1)$ in $f(S^4_2)\x p$ is $K_2$ 
for each $p\in D^2$.

Since $K_2$ and $K_3$ are cobordant, 
there is a compact oriented 3-manifold $P$ in $f(S^4_2)\x I$ 
with the following properties. 
(1)
 $P$ $\cong F\x[0,1]$. 
(2)
$P$ intersects $f(S^4_2)\x\partial I$ transversely. 
$P$ $\cap f(S^4_2)$ in $f(S^4_2)$ is $K_3$.
$P$$\cap [f(S^4_2)\x\{(1,0)\}]$ in $f(S^4_2)\x\{(1,0)\}$ is $K_2$.

When we rotate $f(S^4_2)\x I$ as above, rotate $P$ together.  
Let $Q$ denote what is made from $P$.

Note that $\partial Q$
=$\partial$ 
$\overline{f(S^4_1)\cap (f(S^4_2)\x D^2)})$  
=$\partial$ 
$\overline{f(S^4_1)-(f(S^4_1)\cap (f(S^4_2)\x D^2)})$.   
Then 
$R$=$\overline{f(S^4_1)-f(S^4_1)\cap (f(S^4_2)\x D^2)}$
$\cup Q$ is a 4-sphere embedded in $S^6$.
Put $\widetilde f(S^4_2)=R$. 
This completes the proof.

\head 6
The proof of Theorem 1.5.
\endhead

It is easy to prove that it suffices to prove:
\proclaim{Proposition} 
Let $L=(K_1, K_2)$ be a ($S^2, S^2$)-link and a slice link. 
Then there exists three 4-spheres $S^4_1$, $S^4_2$, and $S^4$ embedded in $S^6$ 
with the following properties. 

(1)$S^4_1\cap S^4_2=\phi$ (2)$(S^4_1\cap S^4, S^4_2\cap S^4)$ in $S^4$ is $L$. 
\endproclaim 

\f{\bf Proof.}  
Let $S^4$$\times$ $D^2$ denote a tubular neighborhood of $S^4$  in $S^6$. 
Put 
$D^2$=$\{(x,y) \vert$ $x^2+y^2\leqq0 \}$.  
Then $S^4$=$S^4\x\{(0,0)\}$.
Put 
$I$=$\{(x,y) \vert$ $0\leqq x\leqq1,$ $y=0\}$.  
We can regard $S^4$$\times$ $D^2$ 
as the result of rotating $S^4\x I$ around the axis $S^4$.  

Since the 2-link $L$ is slice, 
there exists two 3-discs $D^3_1$ and $D^3_2$ in $S^4\x I$ with the following properties. 
(1)
$D^3_1\cap D^3_2$=$\phi$.
(2)
$D^3_i$ intersects $S^4$ transversely. 
$D^3_i$ $\cap S^4$=$\partial D^3_i$.
(3)
($\partial D^3_1$, $\partial D^3_2$) in $S^4$ is the 2-link $L$.

When we rotate $S^4\x I$ as above, 
we rotate   
$ D^3_1$ $\amalg D^3_2$ together. 
This gives 4-spheres $S^4_1$ and $S^4_2$ embedded in $S^6$. 
This completes the proof.

\np

Department of Mathematical Sciences,  University of Tokyo  

  Komaba, Tokyo 153,    Japan

i33992\@m-unix.cc.u-tokyo.ac.jp



\enddocument